\documentclass[12pt,a4paper,english,intoc,bibliography=totoc,index=totoc,BCOR10mm,captions=tableheading,titlepage,fleqn]{scrbook}
\usepackage{lmodern}

\usepackage[T1]{fontenc}
\usepackage[latin9]{inputenc}
\usepackage{babel}
\usepackage{amsmath}
\usepackage{amssymb}
\usepackage[unicode=true,
 bookmarks=true,bookmarksnumbered=true,bookmarksopen=true,bookmarksopenlevel=2,
 breaklinks=false,pdfborder={0 0 0},backref=false,colorlinks=false]
 {hyperref}
\hypersetup{pdftitle={Abstract},
 pdfauthor={Uwe Stöhr},
 pdfpagelayout=OneColumn, pdfnewwindow=true, pdfstartview=XYZ, plainpages=false}
\usepackage{breakurl}

\makeatletter

\usepackage{calc}

\makeatother

\begin{document}
\noindent \begin{center}
\textbf{\Large $\kappa$-Madness and Definability}
\par\end{center}{\Large \par}

\noindent \begin{center}
{\large Haim Horowitz and Saharon Shelah}
\par\end{center}{\large \par}

\noindent \begin{center}
\textbf{\small Abstract}
\par\end{center}{\small \par}

\noindent \begin{center}
{\small Assuming the existence of a supercompact cardinal, we construct
a model where, for some uncountable regular cardinal $\kappa$, there
are no $\Sigma^1_1(\kappa)-\kappa-$mad families.}%
\footnote{{\small Date: May 15, 2018}{\small \par}



Publication 1145 of the second author.%
}
\par\end{center}{\small \par}

\textbf{\large Introduction}{\large \par}

The study of higher analogs of descriptive set theoretic results has
gained considerable attention during the past few years. Recent work
includes new results on regularity properties, definable equivalence
relations and the connections with classification theory (see {[}KLLS{]}
for a survey and a list of relevant open problems).

In this paper we consider the definability of mad families from the
point of view of generalised descriptive set theory. Our basic objects
of study are the following:

\textbf{Definition 1: }a. A family $\mathcal F \subseteq [\kappa]^{\kappa}$
is called $\kappa-$mad if $|A\cap B|<\kappa$ for every distinct
$A,B \in \mathcal F$, and $\mathcal F$ is $\subseteq-$maximal with
respect to this property.

b. We say that $X\subseteq 2^{\kappa}$ is $\Sigma^1_1(\kappa)$ if
there is a tree $T\subseteq \underset{\alpha<\kappa}{\cup} \kappa^{\alpha} \times 2^{\alpha}$
such that $X=\{ \eta \in 2^{\kappa} : $ there is $\nu \in  \kappa^{\kappa}$
such that $(\nu \restriction \alpha, \eta \restriction \alpha) \in T$
for every $\alpha<\kappa \}$.

Following Mathias' classical result that there are no analytic mad
families ({[}Ma{]}), it's natural to investigate the higher analogs
of Mathias' result for a regular uncountable cardinal $\kappa$. It
turns out that under suitable large cardinal assumptions, it's possible
to construct a model where no $\Sigma^1_1(\kappa)-\kappa-$mad families
exist, thus consistently obtaining a higher version of the result
of Mathias.

The main result of the paper is Theorem 10, which will also be stated
here:

\textbf{Main result: }The existence of a regular uncountable cardinal
$\kappa$ such that there are no $\Sigma^1_1(\kappa)-\kappa-$mad
families is consistent relative to a supercompact cardinal.

An important ingredient of the proof is the forcing $\mathbb{Q}_D$
in Definition 3. $\mathbb{Q}_D$ is a $(<\kappa)$-complete forcing
adding a generic subset of $\kappa$ that is almost contained in every
set from the normal ultrafilter $D$ on $\kappa$. We shall prove
that such forcing notions destroy $\Sigma^1_1(\kappa)-\kappa-$mad
families. Using a Laver-indestructible supercompact cardinal, we shall
iterate those forcings to obtain the desired model

The rest of the paper will be devoted to the proof of the above result.

\textbf{\large Proof of the main result}{\large \par}

\textbf{Hypothesis 2: }We fix a measurable cardinal $\kappa$ and
a normal ultrafilter $D$ on $\kappa$.

We shall now define a variant of Mathias forcing:

\textbf{Definition 3: }A. Let $\mathbb{Q}=\mathbb{Q}_{D}^{\kappa}$
be the forcing notion defined as follows:

a. $p\in \mathbb{Q}$ iff $p=(u,A)=(u_p,A_p)$ where $u\in [\kappa]^{<\kappa}$
and $A\in D$.

b. $\leq=\leq_{\mathbb Q}$ is defined as follows: $p\leq q$ iff

1. $u_p \subseteq u_q$.

2. $A_q \subseteq A_p$.

3. $u_q \setminus u_p \subseteq A_p$.

4. $\alpha<\beta$ for every $\alpha \in u_p$ and $\beta \in u_q \setminus u_p$.

B. Let $\underset{\sim}{u}$ be the $\mathbb Q$-name for $\cup \{u_p : p\in \underset{\sim}{G} \}$.

C. $p\leq^{pr} q$ iff $p\leq q$ and $u_p=u_q$.

\textbf{Observation 4: }a. $\mathbb Q$ is $(<\kappa)$-complete.

b. The sequence $(p_i : i<\kappa)$ has an upper bound if the following
conditions holds:

1. $(p_i : i<\kappa)$ is $\leq^{pr}$-increasing.

2. If $i\in \underset{j<i}{\cap}A_j$ and $i>sup(u_{p_0})$ then $j\in [i,\kappa) \rightarrow i\in A_{p_j}$.

\textbf{Proof: }a. By the $\kappa$-completeness of $D$.

b. By the normality of $D$, $(u_{p_0},\underset{i<\kappa}{\Delta}A_{p_i} \setminus u_{p_0})$
is a condition in $\mathbb Q$, it's easy to see that it's the desired
upper bound. $\square$

\textbf{Claim 5:} Suppose that $p\in \mathbb Q$, $sup(u_p) \leq \alpha<\kappa$
and $\underset{\sim}{\tau}$ is a $\mathbb Q$-name of a member of
$V$, then there is $q\in \mathbb Q$ such that:

a. $p\leq^{pr} q$.

b. $A_q \cap (\alpha+1)=A_p \cap (\alpha+1)$.

c. If $v\subseteq \alpha+1$ and there is $r\in \mathbb Q$ forcing
a value to $\underset{\sim}{\tau}$ such that $u_r=v$, then $q^{[v,\alpha]}:=(v,A_q \setminus (\alpha+1)) \in \mathbb Q$
forces the same value to $\underset{\sim}{\tau}$.

\textbf{Proof: }Fix an enumeration $(v_{\beta} : \beta<2^{|\alpha|})$
of $\mathcal P(\alpha+1)$. We shall construct by induction a decreasing
sequence $(A_{\beta} : \beta<2^{|\alpha|})$ of elements of $D$ as
follows:

a. $\beta=0$: Without loss of generality, there is $r\in \mathbb Q$
as in clause (c) for $v_0$. Let $A_0=A_r \cap A_p$.

b. $\beta$ is a limit ordinal: Let $A_{\beta}=\underset{\gamma<\beta}{\cap}A_{\gamma} \in D$
(recall that $2^{|\alpha|}<\kappa$).

c. $\beta=\gamma+1$: Without loss of generality, there is $r\in \mathbb Q$
as in clause (c) for $v_{\beta}$. Let $A_{\beta}=A_{\gamma} \cap A_r$.

Now let $A_q:= ((\underset{\beta<2^{|\alpha|}}{\cap}A_{\beta}) \setminus (\alpha+1)) \cup (A_p \cap (\alpha+1)) \in D$
and $u_q:=u_p$. It's now easy to verify that $q$ is as required.
$\square$

\textbf{Claim 6: }If $p \in \mathbb Q$, $p\Vdash "\underset{\sim}{\tau} \in V^{\kappa}"$
and $sup(u_p) \leq \alpha<\kappa$, then there is $q\in \mathbb Q$
satisfying clause (a) from Claim 5, and in addition: If $i<\kappa$,
$v\subseteq (\alpha+1)$ and there is $r\in \mathbb Q$ forcing a
value to $\underset{\sim}{\tau}(i)$ such that $u_r=v$, then $q^{[v,\alpha]}$
forces the same value to $\underset{\sim}{\tau}(i)$.

\textbf{Proof: }We construct a $\leq^{pr}$ increasing sequence $(p_i : i<\kappa)$
by induction on $i<\kappa$ as follows:

a. i=0: Let $p_0$ be $q$ from the previous claim, where $\underset{\sim}{\tau}(0)$
here stands for $\underset{\sim}{\tau}$ there.

b. $i=j+1$: Similarly, letting $(p_j,\underset{\sim}{\tau}(j))$
here stand for $(p,\underset{\sim}{\tau})$ in Claim 5, let $p_i$
be the corresponding $q$ from Claim 5.

c. $i$ is a limit ordinal: Let $p_i'$ be an upper bound for $(p_j : j<i)$
(see Observation 4). It's easy to see that if the sequence is $\leq^{pr}$-increasing,
then we can get a $\leq^{pr}$-upper bound. Now construct $p_i$ as
in the previous case.

Finally, let $q$ be a $\leq^{pr}$-upper bound for $(p_i : i<\kappa)$
(such $q$ exists by Observation 4(b)). $q$ is obviously as required.
$\square$

\textbf{Claim 7: }If $p \in \mathbb Q$ and $p\Vdash "\underset{\sim}{\tau} \in V^{\kappa}"$,
then there is $q \in \mathbb Q$ that satisfies the conclusion of
Claim 6 for every $\alpha \in [sup(u_p),\kappa)$.

\textbf{Proof: }By Claim 6 and Observation 4(b). $\square$

\textbf{Claim 8: }$(\alpha)$ (A) implies (B) where:

A. a. $\bold B$ is a $\Sigma^1_1(\kappa)$ subset of $[\kappa]^{\kappa}$
and $\Vdash "\underset{\sim}{X} \in \bold B"$.

b. $\chi >2^{\kappa}$, $N\prec (H(\chi),\in)$, $\{\bold{B},D,\underset{\sim}{X}\} \subseteq N$,
$|N|=\kappa$ and $[N]^{<\kappa} \subseteq N$.

c. $\mathbb Q$ is a $(<\kappa)$-complete forcing notion.

d. $\mathbb Q \in N$.

e. $G\subseteq \mathbb Q \restriction N$ is generic over $N$.

B. $\underset{\sim}{X}[G]$ is well defined and belongs to $\bold{B}$.

$(\beta)$ (A) implies (B) where:

A. a. $\bold B$ is a $\Sigma^1_1(\kappa)$ subset of $[\kappa]^{\kappa}$
defined by the tree $T \in V$.

b. $\mathbb Q$ is a $(<\kappa)$-complete forcing notion. 

c. $\bold{B}^{V^{\mathbb Q}}$ is $\kappa$-mad in $V^{\mathbb Q}$.

B. $\bold{B}^V$ is $\kappa$-mad in $V$.

\textbf{Proof: }$(\alpha)$ For $\alpha<\kappa$, let $T_{\alpha}=2^{\alpha} \times \kappa^{\alpha}$,
and for $\alpha<\beta \leq \kappa$ and $(\eta,\nu) \in T_{\beta}$,
let $(\eta,\nu) \restriction \alpha=(\eta \restriction \alpha,\nu \restriction \alpha) \in T_{\alpha}$.
Let $T_*=\underset{\alpha<\kappa}{\cup}T_{\alpha}$, then $T_{\kappa}$
is the set of $\kappa-$branches through $T_*$. There is a subtree
$T\subseteq T_*$ such that $\{\eta : (\eta,\nu) \in lim(T)\}=\bold B$
(where $\eta$ is interpreted as $\{ \alpha: \eta(\alpha)=1\}$),
hence there are $(\underset{\sim}{\eta},\underset{\sim}{\nu})$ such
that $\Vdash "(\underset{\sim}{\eta},\underset{\sim}{\nu}) \in lim(T)$
and $\underset{\sim}{X}=\{ \alpha: \underset{\sim}{\eta}(\alpha)=1\}"$.
Without loss of generality, $\underset{\sim}{\eta},\underset{\sim}{\nu} \in N$.
For each $\alpha<\kappa$, let $I_{\alpha} \in N$ be a dense open
subset of $\mathbb Q$ where $I_{\alpha}=\{ p\in \mathbb Q : p$ forces
a value to $(\underset{\sim}{\eta}, \underset{\sim}{\nu}) \restriction \alpha\}$.
For each $\alpha<\kappa$, choose $p_{\alpha} \in G \cap I_{\alpha}$
and let $(\eta_{\alpha}, \nu_{\alpha}) \in T_{\alpha}$ be the valued
forced by $p_{\alpha}$ for $(\underset{\sim}{\eta},\underset{\sim}{\nu}) \restriction \alpha$.
For every $\alpha<\beta<\kappa$, $p_{\alpha}$ and $p_{\beta}$ are
compatible and hence $\eta_{\alpha} \leq \eta_{\beta}$ and $\nu_{\alpha} \leq \nu_{\beta}$.
Let $(\eta,\nu):=(\underset{\alpha<\kappa}{\cup} \eta_{\alpha},\underset{\alpha<\kappa}{\cup} \nu_{\alpha}) \in lim(T)$,
then $N[G] \models "\underset{\sim}{X}[G]=\{ \alpha: \eta(\alpha)=1\}"$,
hence $\underset{\sim}{X}[G] \in \bold B$. This completes the proof
of $(\alpha)$.

$(\beta)$ Obviously, each element of $\bold{B}^V$ has cardinality
$\kappa$ and $\bold{B}^V$ is a $\kappa$-almost disjoint family.
Let $C\in [\kappa]^{\kappa}$, by assumption (A)(c), $\Vdash_{\mathbb Q} "$there
is $D\in \bold B$ such that $|C\cap D|=\kappa"$. Therefore, for
some $\mathbb Q$-name $\underset{\sim}{\tau}$, $\Vdash_{\mathbb Q} "\underset{\sim}{\tau} \in \bold B$
and $|C\cap \underset{\sim}{\tau}|=\kappa"$. Fix a large enough $\chi$
and $N \prec (H(\chi),\in)$ such that $|N|=\kappa$, $[N]^{<\kappa}$
and $\{\underset{\sim}{\tau},\bold B,C\} \subseteq N$. By the $(<\kappa)$-completeness
of $\mathbb Q$, there is $G\subseteq \mathbb Q \restriction N$ which
is generic over $N$. By part $(\alpha)$ of the claim, $\underset{\sim}{\tau}[G] \in \bold{B}^V$
and $|C\cap \underset{\sim}{\tau}[G]|=\kappa$, hence $\bold{B}^V$
is $\kappa$-mad in $V$. $\square$

\textbf{Claim 9: }There are no $(\mathbb Q,\underset{\sim}{u},D,\bold B)$
such that:

a. $\mathbb Q$ is a $(<\kappa)$-complete forcing notion.

b. $D$ is a normal ultrafilter on $\kappa$.

c. $\Vdash_{\mathbb Q} "\underset{\sim}{u} \in [\kappa]^{\kappa}$
and $\underset{\sim}{u} \subseteq^* A$ for every $A\in D"$.

d. $\bold B \in V$ is a $\Sigma^1_1(\kappa)$ subset of $[\kappa]^{\kappa}$.

e. $\bold B^V$ is $\kappa$-mad in $V$.

f. $\bold B^{V^{\mathbb Q}}$ is $\kappa$-mad in $V^{\mathbb Q}$.

\textbf{Proof: }Suppose towards contradiction that there are $(\mathbb Q,\underset{\sim}{u},D,\bold B)$
as above. Hence $\bold B$ is a $\Sigma^1_1(\kappa)$-$\kappa$-mad
family in $V$. Fix a sequence $(A_i^* : i<\kappa) \in V$ of pairwise
distinct members of $\bold B$. Let $F: \kappa \times \kappa \rightarrow \kappa$
be the function defined as $F(i,\alpha):=$the $\alpha$th member
of $A_i^* \setminus \underset{j<i}{\cup} A_j^* \in [\kappa]^{\kappa}$
(recalling that $\kappa$ is regular and $\bold B$ is $\kappa$-almost
disjoint).

Now define the following $\mathbb Q$-names:

1. $\underset{\sim}{\alpha_i}$ is $min\{ \underset{\sim}{u} \setminus (i+1)\}$.

2. $\underset{\sim}{\beta_i}$ is $F(i,\underset{\sim}{\alpha_i})$.

3. $\underset{\sim}{v}=\{ \underset{\sim}{\beta_i} : i\in \underset{\sim}{u}$
satisfies that $otp(i\cap \underset{\sim}{u})$ is even $\}$.

Let $E$ be the ultrafilter on $\kappa$ generated by the sets $\{ \{F(i,\alpha) : i<\alpha$
are from $A \} : A\in D\}$. By Rowbottom's theorem, for every $A\in D$
and $X \subseteq \kappa$, if $f_X : [A]^2 \rightarrow \{0,1\}$ is
defined by $f_X(i,\alpha)=0$ iff $F(i,\alpha) \in X$, then there
exists a monochromatic $B\subseteq A$ such that $B\in D$. It follows
that $E$ is indeed an ultrafilter. As $F$ is injective, each set
in $E$ has cardinality $\kappa$. By the $\kappa$-completeness of
$D$, $E$ is also $\kappa$-complete.

Subclaim 1: $E \cap \bold B=\emptyset$.

Proof: Let $C\in \bold B$.

Case I: $C=A_j^*$ for some $j<\kappa$. Let $A\in D$ such that $min(A)>j$,
then by the definition off $F$, $\{F(i,\alpha) : i<\alpha$ are from
$A \} \cap A_j^*=\emptyset$. It follows that $C \notin E$.

Case II: $C \in \bold B \setminus \{A_i^* : i<\kappa\}$. In this
case, define $f:\kappa \rightarrow \kappa$ by $f(i)=sup(A_i^* \cap C)+i+1$
and let $H=\{ \delta<\kappa : \delta$ is a limit ordinal such that
$f(i)<\delta$ for all $i<\delta \}$. $H\subseteq \kappa$ is a club,
hence $H\in D$ and $H^*:=\{F(i,\alpha) : i<\alpha$ are from $H\} \in E$.
Suppose that $F(i,\alpha) \in H^*$, if $F(i,\alpha) \in C$ then
$\alpha \leq F(i,\alpha)<f(i)<\alpha$, a contradiction. It follows
that $C\notin E$.

This proves the subclaim. We shall now return to the proof of the
main claim. Suppose towards contradiction that $\bold{B}^{V^{\mathbb Q}}$
is $\kappa-$mad in $V^{\mathbb Q}$. As $\Vdash_{\mathbb Q} "\underset{\sim}{v} \in [\kappa]^{\kappa}$,
there is a $\mathbb Q$-name $\underset{\sim}{\tau}$ of a member
of $\bold{B}^{V^{\mathbb Q}}$ such that $\Vdash_{\mathbb Q} "|\underset{\sim}{v} \cap \underset{\sim}{\tau}|=\kappa"$.
For every $p\in \mathbb Q$, let $B_p^+=\{ \alpha<\kappa : p \nVdash "\alpha \notin \underset{\sim}{\tau}"\}$.

Subclaim 2: $B_p^+ \in E$.

Proof: Suppose towards contradiction that $B_p^+ \notin E$, then
there is some $C_p \in D$ such that $B_p^+ \cap \{F(i,\alpha) : i<\alpha$
are from $C_p\}=\emptyset$. Therefore, if $i<\alpha$ are from $C_p$
then $p\Vdash "F(i,\alpha) \notin \underset{\sim}{\tau}"$. Recalling
that $\Vdash_{\mathbb Q} "\underset{\sim}{u} \subseteq^* C_p"$, it
follows that $p\Vdash "\underset{\sim}{\alpha_i} \in C_p$ for $i$
large enough$"$, and also $p\Vdash "$for $i$ large enough, $i \in \underset{\sim}{u} \rightarrow i\in C_p"$.
Therefore, $p\Vdash_{\mathbb Q} "\underset{\sim}{\beta_i}=F(i,\underset{\sim}{\alpha_i}) \notin \underset{\sim}{\tau}$
for every large enough $i\in \underset{\sim}{u}"$. Recalling the
definition of $\underset{\sim}{v}$, it follows that $p\Vdash "|\underset{\sim}{v} \cap \underset{\sim}{\tau}|<\kappa"$,
contradicting the choice of $\underset{\sim}{\tau}$. It follows that
$B_p^+ \in E$, which completes the proof of Subclaim 2.

For every $p\in \mathbb Q$, let $B_p^-=\{ \alpha<\kappa : p\nVdash "\alpha \in \underset{\sim}{\tau}"\}$.

Subclaim 3: $B_p^- \in E$.

Proof: Suppose not, then $B_*:=\kappa \setminus B_p^- \in E$ (hence
$B_* \in [\kappa]^{\kappa}$) and $p\Vdash "B_* \subseteq \underset{\sim}{\tau}"$.
By the $\kappa$-madness of $\bold B$, there is $C\in \bold B$ (in
$V$) such that $|C\cap B_*|=\kappa$. As $p\Vdash "B_* \cap C \subseteq \underset{\sim}{\tau},$
$\underset{\sim}{\tau} \in \bold B$ and $\bold B$ is $\kappa$-mad$"$,
it follows that $p\Vdash "\underset{\sim}{\tau}=C"$. We shall derive
a contradiction by showing that $\Vdash_{\mathbb Q} "|\underset{\sim}{\nu} \cap C|<\kappa"$:
Choose $i_*$ such that $C\neq A_i^*$ for every $i\in [i_*,\kappa)$.
It follows that $|C\cap A_i^*|<\kappa$ for every $i\in [i_*,\kappa)$.
Now repeat the argument of Case II in the proof of Subclaim 1 and
choose $f$, $H$ and $H^*$ as there. As $H\in D$, $\Vdash_{\mathbb Q} "$for
large enough $i$, $i\in \underset{\sim}{u} \rightarrow i, \underset{\sim}{\alpha_i}\in H"$.
Repeating the same argument as in Subclaim 1, $\Vdash_{\mathbb Q} "$for
large enough $i\in \underset{\sim}{u}$, $\underset{\sim}{\beta_i}=F(i,\underset{\sim}{\alpha_i}) \in H^*$,
hence $\underset{\sim}{\beta_i} \notin C"$. It follows that $\Vdash_{\mathbb Q} "|\underset{\sim}{v} \cap C|<\kappa"$,
leading to a contradiction. This completes the proof of Subclaim 3.

Observation 4: A. Given $p_1,p_2 \in \mathbb Q$ and $\alpha<\kappa$,
there exist $(q_1,q_2,\beta)$ such that:

a. $p_l \leq_{\mathbb Q} q_l$ $(l=1,2)$.

b. $\beta \in [\alpha,\kappa)$.

c. $p_1 \Vdash "\beta \in \underset{\sim}{\tau}"$.

d. $p_2 \Vdash "\beta \notin \underset{\sim}{\tau}"$.

B. As in (A), with (d) replaced by the following:

d'. $p_2 \Vdash "\beta \in \underset{\sim}{\tau}"$.

Proof: By the previous subclaims, $B_{p_1}^+ \cap B_{p_2}^-, B_{p_1}^+ \cap B_{p_2}^+ \in E$,
hence there exist $\beta \in (B_{p_1}^+ \cap B_{p_2}^-) \setminus \alpha$
and $\gamma \in (B_{p_1}^+ \cap B_{p_2}^+) \setminus \alpha$. By
the definitions of $B_p^{+/-}$, there exist $q_1 \geq p_1$ and $q_2 \geq p_1$
such that $(q_1,q_2,\beta)$ are as required, and similarly for $\gamma$
and (B). This proves the observation.

Let $\chi=(2^{\kappa})^+$ and $N\prec (H(\chi),\in)$ such that $|N|=\kappa$,
$N^{<\kappa} \subseteq N$, $\kappa \subseteq N$ and $\underset{\sim}{\tau},D,\bold B\in N$.
Let $(I_i : i<\kappa)$ list the dense open subsets of $\mathbb Q$
from $N$. We shall now choose $(p_i^1,p_i^2,\beta_i)$ by induction
on $i<\kappa$ such that:

a. $p_i^1,p_i^2 \in \mathbb Q \cap N$ and $\beta_i \in N$.

b. $i<j \rightarrow p_i^l \leq_{\mathbb Q} p_j^l$ $(l=1,2)$.

c. If $i=4j+1$ then $p_i^1,p_i^2 \in I_j$.

d. $\beta_i \in \kappa \setminus \underset{j<i}{\cup} (\beta_j+1)$.

e. If $i=4j+2$ then $p_i^1 \Vdash "\beta_{4j+2} \in \underset{\sim}{\tau}"$
and $p_i^2 \Vdash "\beta_{4j+2} \in \underset{\sim}{\tau}"$.

f. If $i=4j+3$ then $p_i^1 \Vdash "\beta_{4j+3} \in \underset{\sim}{\tau}"$
and $p_i^2 \Vdash "\beta_{4j+3} \notin \underset{\sim}{\tau}"$.

g. If $i=4j+4$ then $p_i^1 \Vdash "\beta_{4j+4} \notin \underset{\sim}{\tau}"$
and $p_i^2 \Vdash "\beta_{4j+4} \in \underset{\sim}{\tau}"$.

Observation 5: It is possible to choose $(p_i^1,p_i^2,\beta_i)$ as
above for each $i<\kappa$.

Proof:

Case I: $i=0$. This is trivial.

Case II: $i$ is a limit ordinal: As $N^{<\kappa} \subseteq N$ and
$(p_j^l : j<i), (\beta_j :j<i) \in N$, we can find $p_i^1$ and $p_i^2$
using the $(<\kappa)$-completeness of $\mathbb Q$ and elementarity.
As $\kappa$ is regular, there is no problem to choose $\beta_i$.

Case III: $i=4j+1$: As $p_j^1,p_j^2,I_j \in N$, by elementarity
there exist $p_i^1$ and $p_i^2$ as required.

Case IV: $i=4j+2$: Use Observation 4(B).

Case V: $i=4j+3$: Use Observation 4(A).

Case VI: $i=4j+4$: Use Observation 4(A), with $(p_i^2,p_i^1)$ here
standing for $(p_1,p_2)$ there.

Finally, let $G_l=\{q\in \mathbb Q \cap N : q\leq_{\mathbb Q}p_i^l$
for some $i<\kappa \}$ $(l=1,2)$, then $G_l \subseteq \mathbb Q \cap N$
is generic over $N$. By Claim $8(\alpha)$, $C_l:=\underset{\sim}{\tau}[G_l] \in \bold B$.
By the choice of $(p_i^1,p_i^2,\beta_i)$, $\{\beta_{4i+2} : i<\kappa\} \subseteq C_1 \cap C_2$,
hence $C_1 \cap C_2 \in [\kappa]^{\kappa}$. Similarly, $|\{ \beta_{4i+3} : i<\kappa\}|=\kappa$
and $\{ \beta_{4i+3} : i<\kappa\} \subseteq C_1 \setminus C_2$, hence
$C_1 \neq C_2$. This contradicts the $\kappa$-madness of $\bold B$
in $V$, which completes the proof of Claim 9. $\square$

\textbf{Theorem 10: }If $\kappa$ is a Laver-indestructible supercompact
cardinal then there is a generic extension where $\kappa$ is supercompact
and there are no $\Sigma^1_1(\kappa)$-$\kappa$-mad families.

\textbf{Proof: }We recall the following strong version of $\kappa^+-c.c.$
(see e.g. {[}Sh:80{]} and {[}Sh:1036{]}): A forcing $\mathbb Q$ satisfies
$*^1_{\mathbb \kappa,\mathbb Q}$ if:

a. $\mathbb Q$ is $(<\kappa)$-complete.

b. If $\{p_{\alpha} : \alpha<\kappa^+\} \subseteq \mathbb Q$, then
for some club $E\subseteq \kappa^+$ and pressing down function $f$
on $E$ we have $(\delta_1,\delta_2 \in E \wedge f(\delta_1)=f(\delta_2)) \rightarrow p_{\delta_1},p_{\delta_2}$
are compatible.

c. Every two compatible conditions in $\mathbb Q$ have a least upper
bound.

Obviously, $*^1_{\kappa,\mathbb Q}$ implies $\kappa^+-c.c.$. By
{[}Sh:80{]}, $*^1_{\kappa,\mathbb Q}$ is preserved under $(<\kappa)$-support
iterations.

It's easy to verify that $\mathbb Q= \mathbb{Q}_D$ satisfies $*^1_{\kappa,\mathbb Q}$
when $D$ is a normal ultrafilter on $\kappa$ (e.g. fix a bijection
$g: [\kappa]^{<\kappa} \rightarrow \kappa$, and for every $\{p_{\alpha} : \alpha<\kappa^+\}$,
let $E=(\kappa,\kappa^+)$ and let $f:E\rightarrow \kappa^+$ be defined
by $f(\alpha)=g(u_{\alpha})$ where $p_{\alpha}=(u_{\alpha},A_{\alpha})$)

Let $(\mathbb{P}_{\alpha},\underset{\sim}{\mathbb{Q}_{\beta}} : \alpha \leq \delta, \beta<\delta)$
be a $(<\kappa)$-support iteration such that:

a. $cf(\delta)>\kappa$.

b. Each $\underset{\sim}{\mathbb{Q}_{\beta}}$ is $*^{1}_{\kappa,\mathbb{Q_{\underset{\sim}{\beta}}}}$.

c. $\delta=sup\{ \alpha<\delta : $ in $V^{\mathbb{P}_{\alpha}}$,
$\underset{\sim}{\mathbb{Q}_{\alpha}}=\mathbb{Q}_{\underset{\sim}{D_{\alpha}}}$
where $\underset{\sim}{D_{\alpha}}$ is a $\mathbb{P}_{\alpha}$-name
of a normal ultrafilter on $\kappa\}$.

As $\kappa$ is a Laver indestructible supercompact cardinal, there
is an iteration as above. Suppose towards contradiction that there
is a $\Sigma^1_1(\kappa)-\kappa-$mad family $\bold B$ in $V^{\mathbb{P}_{\delta}}$.
$\bold B=\{\eta : (\eta,\nu) \in lim(T)\}$ for a suitable tree $T$.
By the fact that $cf(\delta)>\kappa$ and $\mathbb{P}_{\delta}$ is
$\kappa^+-c.c.$, it follows that $T\in V^{\mathbb{P}_{\beta}}$ for
some $\beta<\delta$. Let $\gamma \in [\beta,\delta)$ such that $\underset{\sim}{\mathbb{Q}_{\gamma}}=\mathbb{Q}_{\underset{\sim}{D_{\gamma}}}$
where $\underset{\sim}{D_{\gamma}}$ is a $\mathbb{P}_{\gamma}$-name
of a normal ultrafilter on $\kappa$. By Claim $8(\beta)$, $\bold{B}^{V^{\mathbb{P}_{\gamma}}}$
is $\kappa-$mad in $V^{\mathbb{P}_{\gamma}}$.

Applying Claim 9 to $V_1=V^{\mathbb{P}_{\gamma}}$, $\mathbb{Q}=\mathbb{P}_{\delta}/\mathbb{P}_{\gamma}$
and $D=\underset{\sim}{D}_{\gamma}$, it follows that $\bold B$ is
not $\kappa-$mad in $V^{\mathbb{P}_{\delta}}$, a contradiction.
It follows that there are no $\Sigma^1_1(\kappa)-\kappa-$mad families
in $V^{\mathbb{P}_{\delta}}$. $\square$

\textbf{\large Open problems}{\large \par}

We conclude by listing some of the open problems following from our
work:

Following the main result of the paper, one may ask whether it's possible
to get an implication instead of just consistency:

\textbf{Question 1: }Suppose that $\kappa$ is supercompact, is there
a $\Sigma^1_1(\kappa)-\kappa$-mad family?

\textbf{Question 2: }What is the consistency strength of $ZFC+"$for
some uncountable regular cardinal $\kappa$, there are no $\Sigma^1_1(\kappa)-\kappa-$mad
families$"$?

It's known by {[}Ma{]}, {[}To{]} and {[}HwSh:1090{]} that $ZF+DC+"$there
are no mad families$"$ is consistent ({[}To{]} shows that it holds
in Solovay's model while in {[}HwSh:1090{]} we obtain a consistency
result relative to $ZFC$).

\textbf{Question 3: }a. What's the consistency strength of $ZF+DC+$there
exists a regular uncountable cardinal $\kappa$ such that there are
no $\kappa-$mad families$"$?

b. Suppose that $\kappa>\aleph_0$ is regular, does $DC_{\kappa}$
imply the existence of a $\kappa-$mad family?

It's known by {[}HwSh:1089{]} and {[}HwSh:1095{]} that Borel maximal
eventually different families and maxima cofinitary groups exist,
therefore it's natural to investigate the $\kappa$-version of those
results:

\textbf{Question 4: }a. Does $ZFC$ imply that there are $\kappa-$Borel
$\kappa-$maximal eventually different families for every (or at least
for some) regular uncountable cardinal $\kappa$?

b. Similarly, replacing regular uncountable cardinals by successor
cardinals, inaccessible non-Mahlo cardinals, etc.

\textbf{Question 5: }Does $ZFC$ imply that there are $\kappa-$Borel
$\kappa-$maximal cofinitary groups for every (or at least for some)
regular uncountable cardinal $\kappa$?

b. Similarly, replacing regular uncountable cardinals by successor
cardinals, inaccessible non-Mahlo cardinals, etc.

\textbf{\large References}{\large \par}

{[}HwSh:1089{]} Haim Horowitz and Saharon Shelah, A Borel maximal
eventually different family, arXiv:1605.07123

{[}HwSh:1090{]} Haim Horowitz and Saharon Shelah, Can you take Toernquist's
inaccessible away?, arXiv:1605.02419

{[}HwSh:1095{]} Haim Horowitz and Saharon Shelah, A Borel maximal
cofinitary group, arXiv:1610.01344

{[}KLLS{]} Yurii Khomskii, Giorgio Laguzzi, Benedikt Loewe and Ilya
Sharankou, Questions on generalised Baire spaces, Math. Log. Quart.
\textbf{62}, No. 4-5, 439-456 (2016)

{[}Ma{]} A. R. D. Mathias, Happy families, Ann. Math. Logic \textbf{12
}(1977), no. 1, 59-111

{[}Sh:80{]} Saharon Shelah, A weak generalization of MA to higher
cardinals, Israel J. Math. 30 (1978) 297-306

{[}Sh:1036{]} Saharon Shelah, Forcing axioms for $\lambda-$complete
$\mu^+-c.c.$, arXiv:1310.4042

{[}To{]} Asger Toernquist, Definability and almost disjoint families,
arXiv:1503.07577

$\\$

(Haim Horowitz) Department of Mathematics

University of Toronto

Bahen Centre, 40 St. George St., Room 6290

Toronto, Ontario, Canada M5S 2E4

E-mail address: haim@math.toronto.edu

$\\$

(Saharon Shelah) Einstein Institute of Mathematics

Edmond J. Safra Campus,

The Hebrew University of Jerusalem.

Givat Ram, Jerusalem 91904, Israel.

Department of Mathematics

Hill Center - Busch Campus,

Rutgers, The State University of New Jersey.

110 Frelinghuysen Road, Piscataway, NJ 08854-8019 USA

E-mail address: shelah@math.huji.ac.il
\end{document}